\documentclass[a4paper, 10pt, american]{amsart}

\usepackage{times,latexsym,amssymb}
\usepackage{amsmath,amsthm,bm}
\usepackage{color}
\usepackage[colorlinks,pdfpagelabels,pdfstartview = FitH,bookmarksopen= true,bookmarksnumbered = true,linkcolor = blue,plainpages = false,hypertexnames = false,citecolor = red,pagebackref=false]{hyperref}

\usepackage{amsbsy}
\usepackage{amstext}
\usepackage{amssymb}
\usepackage{esint}
\usepackage{stmaryrd}
\usepackage[pdftex]{graphicx}
\usepackage{floatflt}
\usepackage{appendix}
\usepackage{mathrsfs}
\usepackage{enumerate}

\allowdisplaybreaks

\newtheorem{theorem}{Theorem}[section]
\newtheorem{proposition}{Proposition}[section]
\newtheorem{lemma}{Lemma}[section]
\newtheorem{corollary}{Corollary}[section]
\newtheorem{remark}{Remark}[section]

\numberwithin{equation}{section}
\numberwithin{theorem}{section}
\numberwithin{proposition}{section}
\numberwithin{lemma}{section}
\numberwithin{remark}{section}
\newcommand{\al}{\alpha}

\newcommand{\gm}{\gamma}
\newcommand{\dl}{\delta}

\newcommand{\lm}{\lambda}

\newcommand{\varep}{\varepsilon}

\newcommand{\vp}{\varphi}
\newcommand{\sig}{\sigma}

\newcommand{\om}{\omega}
\newcommand{\Om}{\Omega}

\newcommand{\df}[1]{\buildrel\mbox{\small def}\over{#1}}

\newcommand{\nn}{\mathbb{N}}

\newcommand{\rr}{\mathbb{R}}
\newcommand{\rn}{\rr^N}

\newcommand{\ov}[1]{\overline{#1}}

\newcommand{\bbox}{\vrule height.6em width.6em 
depth0em} 

\newcommand{\dvg}{\operatorname{div}}

\newcommand{\loc}{\operatorname{loc}}

\newcommand{\noi}{\noindent}

\newcommand{\dsty}{\displaystyle}

\def\XXiint#1#2#3{\setbox0=\hbox{$#1{#2#3}{\iint}$}
    \vcenter{\hbox{$#2#3$}}\kern-0.5\wd0}

\numberwithin{equation}{section}

\newtheorem{definition}[theorem]{Definition}

\newcommand{\beq}{\begin{equation}}
\newcommand{\eeq}{\end{equation}}

\def\Xint#1{\mathchoice
    {\XXint\displaystyle\textstyle{#1}}%
    {\XXint\textstyle\scriptstyle{#1}}%
    {\XXint\scriptstyle\scriptscriptstyle{#1}}%
    {\XXint\scriptscriptstyle\scriptscriptstyle{#1}}%
    \!\int}
\def\XXint#1#2#3{\setbox0=\hbox{$#1{#2#3}{\int}$}
    \vcenter{\hbox{$#2#3$}}\kern-0.5\wd0}
\def\bint{\Xint-}
\def\bint{\Xint-}
\def\dashint{\Xint{\raise4pt\hbox to7pt{\hrulefill}}}
\def\dashiint{\bint\kern-0.15cm\bint}

\date{}

\subjclass[2020]{Primary 31B05, 35J25; Secondary 31B15, 35J92}
\keywords{De Giorgi Classes, Rigidity, Sub-Linear Growth, Phragm\'en--Lindel\"of Theorem, Weak Harnack Inequality}

\begin{document}
\title{Phragm\'en-Lindel\"of-type theorems for functions in Homogeneous De Giorgi Classes}
\date{\today}

\author[S. Ciani]{Simone Ciani}
\address{Simone Ciani\\
Dipartimento di Matematica,
Universit\`a di Bologna\\ 
Piazza di Porta San Donato 5, 40126 Bologna, Italy}
\email{simone.ciani3@unibo.it}

\author[U. Gianazza]{Ugo Gianazza}
\address{Ugo Gianazza\\
Dipartimento di Matematica ``F. Casorati",
Universit\`a di Pavia\\ 
via Ferrata 5, 27100 Pavia, Italy}
\email{ugogia04@unipv.it}

\author[Zheng Li]{Zheng Li}
\address{Zheng Li\\
Dipartimento di Matematica ``F. Casorati",
Universit\`a di Pavia\\ 
via Ferrata 5, 27100 Pavia, Italy\\
and
School of Mathematics, Jilin University,
2699 Qianjin Street, Changchun, 130012, P.R. China}
\email{lizheng22@mails.jlu.edu.cn}

\begin{abstract}
We study Phragm\'en-Lindel\"of-type theorems for functions $u$ in homogeneous De Giorgi classes, and we show that the maximum modulus $\mu_+(r)$ of $u$ has a power-like growth of order $\al\in(0,1)$ when $r\to\infty$. By proper counterexamples, we show that in general we cannot expect $\al$ to be $1$.
\end{abstract}

\maketitle

\tableofcontents

\section{Introduction and Main Results}
Let $\Om\subset\rn$ be an open set, and for $x_o\in\rn$ let $B_\rho(x_o)$ denote the ball of radius $\rho$ centered at $x_o$. We have the following well-known definitions.
\begin{definition}\label{DGClass}
Given an open subset $\Om\subseteq\rn$, we say that a function $u \in W^{1,p}_{\loc}(\Om)$ belongs to the \emph{De Giorgi class} $[DG]^{\pm}_p(\Om,\hat\gm)$ for $1<p\le N$ if there exists a constant $\hat{\gamma}>1$ such that for all $\sigma \in (0,1)$, $\rho>0$, $x_o \in \Om$, $k \in \rr$, and $B_\rho(x_o)\subset\Om$,  the following estimate is satisfied:
\begin{equation}\label{DG}
    \int_{B_{\sigma \rho}(x_o)} |D  [(u-k)_{\pm}]|^p \, dx \leq \frac{\hat{\gamma}}{[(1-\sigma) \rho]^p} \int_{B_{\rho}(x_o)} |(u-k)_{\pm}|^p\, dx\,.
\end{equation} 
%
Furthermore, we let
$$[DG]_p(\Om,\hat\gm)\equiv[DG]^{+}_p(\Om,\hat\gm)\cap[DG]^{-}_p(\Om,\hat\gm).$$
\end{definition}
\begin{definition}\label{DGClass-bd}
Let $h\in W^{1,p}(\rn)\cap C(\rn)$. Given an open subset $\Om\subset\rn$, the De Giorgi classes $[DG]^{+}_p(\bar\Om,h,\hat\gm)$ in the closure of $\Om$ are defined to consist of functions $u\in [DG]^{+}_p(\Om,\hat\gm)$ such that $(u-h)_+\in W^{1,p}_o(\Om)$ and satisfying in addition inequalities of the type
\[
\forall y\in\partial\Om,\quad \forall B_\rho(y),\quad\forall\,
\sig\in(0,1),\quad\forall  k\geq\sup_{B_\rho(y)\cap\partial\Om}h,
\]
\begin{equation*}\label{DG-boundary}
    \int_{B_{\sigma \rho}(y)\cap\Om} |D  [(u-k)_+]|^p \, dx \leq \frac{\hat{\gamma}}{[(1-\sigma) \rho]^p} \int_{B_{\rho}(y)\cap\Om} |(u-k)_{+}|^p\, dx\,.
\end{equation*}
The various constants are as in \eqref{DG}. The classes $[DG]^{-}_p(\bar\Om,h,\hat\gm)$ are defined analogously by replacing $(u-k)_+$ with $(u-k)_-$ provided that $\dsty k\leq\inf_{B_\rho(y)\cap\partial\Om}h$.
Furthermore, we let
$$[DG]_p(\bar\Om,h,\hat\gm)\equiv[DG]^{+}_p(\bar\Om,h,\hat\gm)\cap[DG]^{-}_p(\bar\Om,h,\hat\gm).$$
\end{definition}
\begin{remark}
\upshape In Definitions~\ref{DGClass}--\ref{DGClass-bd} we referred to balls $B_\rho$; however, if for $y\in\rn$ we 
let $K_\rho(y)$ denote a cube of edge $2\rho$ 
centered at $y$, equivalent definitions can be given in terms of $K_\rho(y)$, instead of balls. Moreover, when the center of a ball or of a cube coincides with the origin, we simply write $B_\rho$ or $K_\rho$.
\end{remark}
A fundamental result by De Giorgi \cite{degiorgi} (see also the work of O.A. Ladyzenskaja \&  N.N. Ural'ceva \cite{ladis}) shows that functions in $[DG]^\pm_p(\Om,\hat\gm)$ are locally bounded and locally H\"older continuous. Moreover, whenever they are non-negative, they satisfy the Harnack inequality \cite{DiBenedetto-Trudy}.

It is important to recall that the De Giorgi classes include local, weak solutions to elliptic equations with bounded and measurable coefficients, subject to some upper and lower ellipticity conditions, but also minima or even $Q$-minima of rather general functionals, even if they do not admit an Euler equation.

Hence, the interest in these classes stems from the fact that their members have very different origins, but they all enjoy properties which are typically associated to solutions to elliptic partial differential equations.

The purpose of this work is precisely in this same direction: we study Phragm\'en-Lindel\"of-type properties satisfied by functions in De Giorgi classes.

We denote half spaces by 
\[
\rn_+= \rr^{N-1} \times \{x_N>0\}, \qquad \text{and} \qquad \rn_-= \rr^{N-1} \times \{x_N< 0\}.
\]

In the sequel, we are interested in taking either $\Om=\rn_+$ and $u\in[DG]^{+}_p(\bar\Om,h,\hat\gm)$ such that $h=0$ on $\partial\rn_+$, or $\Om=\rn$ and $u\in[DG]^{+}_p(\rn,\hat\gm)$ vanishing in a proper subset of $\rn$. In both cases, since the boundary value of $u$ is zero, we need to take $k\ge0$; therefore, it is not restrictive to directly assume $u\ge0$.

Our main results are the following theorems.
\begin{theorem}\label{PL} 
Take $1<p\le N$ and let $u \in [DG]^{+}_p(\rn_+,\hat\gm)$; assume that $u$ is non-negative, $u\not\equiv 0$ in $\rn_+$, and such that 
\begin{equation} \label{boundary-condition}
    \lim_{\rn_+ \ni \, y \rightarrow x \, \in \partial \rn_+} u\, = 0.
\end{equation} 
Then, there exists $0<\al_1< 1$ depending on $N$, $p$, and $\hat\gm$, such that  
\begin{equation*} 
0<\liminf_{r\uparrow \infty}\,  \frac{\mu_+(r)}{r^{\al_1}}, \quad \text{ where} \quad \mu_+(r)= \sup_{B_r \cap \rn_+} u \,.
\end{equation*}
\end{theorem}
\begin{remark}
Condition \eqref{boundary-condition} is to be understood in the sense of Definition~\ref{DGClass-bd} with $h\equiv0$.
\end{remark}
\begin{theorem}\label{PL-2} 
Take $1<p\le N$, let $u \in [DG]^{+}_p(\rn,\hat\gm)$, assume that $u$ is non-negative, $u\not\equiv0$ in $\rn$, and $u$ vanishes in 
${\mathbb H}^{M}$, the $M$-dimensional hyperplane through the origin, that is
\[ 
{\mathbb H}^{M}\df=\{(x_1,\dots,x_M,0,0,\dots,0)\}\quad\text{ with }\, 1\le M\le N-1.
\] 
Then, if $N-M<p\le N$, there exists $0<\al_2<1$ depending on $N$, $p$, $M$, and $\hat\gm$, such that 
\begin{equation*} 
0<\liminf_{r\uparrow \infty}\,  \frac{\mu_+(r)}{r^{\al_2}}, \quad \text{ where} \quad \mu_+(r)= \sup_{B_r \cap \rn} u \,.
\end{equation*}
\end{theorem}
\begin{remark}
{\normalfont Since functions $u\in W^{1,p}(\Om)$ are defined up to set of \emph{$p$-capacity} zero (see below \S~\ref{SS:cap1} for the definition of this notion), assuming that $u$ vanishes on ${\mathbb H}^M$ is well-defined.}
\end{remark}  

\subsection{Novelty and Significance}
It is quite difficult to give an exhaustive and complete overview of the state of the art concerning Phragm\'en–Lindel\"of-type theorems. In the following, we limit ourselves to just few results, which are more directly connected to our work. 

The first and classical formulation (see \cite{PH}) deals with the growth of a subharmonic function $u$ in a halfplane $\rr^2_+$; the theorem states that if $u$ is non-positive on the boundary of $\rr^2_+$ and admits a positive value in it, then
\[
\liminf_{r\uparrow\infty}\frac{\mu_+(r)}r>0,
\]
where $\displaystyle\mu_+(r)=\sup_{(x,y)\in B_r\cap \rr^2_+}|u(x,y)|$.
\vskip0.2cm 
\noindent Even though the original result was stated for functions of a complex variable $z=x+iy$, nevertheless, it is essentially a theorem on harmonic functions of two variables. It was later generalized to harmonic functions of $N$ variables by Ahlfors (see \cite[Theorem~8]{Ahlfors}, and also the comment immediately below the statement). In \cite{Granlund-1985} Granlund, Lindqvist and Martio studied the conformally invariant case, i.e., they considered extremals of the variational integral
\[
\int_\Om F(x,Du)dx,
\]
where $F(x,\xi)\approx|\xi|^N$, $\xi \in \rn$ (the precise assumptions are collected in \S~2.1 of their paper). 
\vskip0.2cm 
\noindent Thus, the plane harmonic case is included in their study, but the classical harmonic case in $\rn$ with $N\ge3$ is not. In this wider context, they obtained a further, proper generalization.

A Phragm\'en-Lindel\"of-type theorem for non-negative sub-solutions to the parabolic $p$-Laplacian (so-called $p$-subharmonic functions) is given in \cite[Theorem~4.6]{Lindqvist1985}; in the same paper it is also studied the behavior of $N$-subharmonic functions which vanish on $\mathbb H^q$ for $1\le q\le N-1$ \cite[Theorem~4.8]{Lindqvist1985}. The proofs are based on the comparison principle, and on explicit estimates of the $N$-harmonic measure. The sharpness of these results is discussed with the aid of straightforward counterexamples.

Although explicit computations of the $p$-harmonic measure for a general $p>1$ are not as easy as for the $p=N$ case, nevertheless, it turns out that proper estimates from below and from above suffice, and this allowed Lundstr\"om to extend \cite[Theorem~4.8]{Lindqvist1985} to a general $p$ (see \cite[Corollary~4.3]{Lund}, and obtain the corresponding Phragm\'en-Lindel\"of-type theorem. 

There is another important fact which is worth pointing out in this context, even though we will not touch upon it in the following. Once more, we do not pretend to give a full account of the state of the art.

\noindent If one directly assumes that $u$ is a non-negative solution to $\Delta_p u = 0$ in $\{x_N> 0\}$ which continuously vanishes on the flat boundary $\{x_N = 0\}$, then Phragm\'en-Lindel\"of \emph{qualitative} results can be \emph{quantified}, and one can show that, modulo a proper normalization, $\displaystyle u(x) = x_N$ in $\{x_N\ge0\}$. Such a result was probably stated for the first time in a paper due to Loomis \& Widder \cite{Loomis} about two-variable harmonic functions. Then, it was extended to positive harmonic functions in $N$ variables by Rudin \cite[Lemma~1]{Rudin}, whose proof is based on the same ideas of \cite{Loomis}. Related results were obtained by Gilbarg and Serrin in the plane for solutions to linear elliptic equations in non-divergence form with smooth coefficients (see \cite{Gilbarg,Serrin}).

The corresponding extension to non-negative solutions of the $p$-Laplacian was studied, for example, in \cite{KSZ}. An interesting discussion of further extensions is given in the introduction of \cite{Moreira}.

Coming back to the qualitative formulation, we were inspired by the results in \cite{Lindqvist1985,Lund}. Due to the wide generality of functions in the De Giorgi classes, there is no chance of using the comparison principle in our case (indeed, whether any form of comparison principle holds is a major open problem for functions in $[DG]^\pm_p(\Om,\hat\gm)$, even though it is generally believed that such a result does not hold). Hence, we had to develop a different approach, based on a proper use of the weak Harnack inequality (see Lemma~\ref{lemma-WH} and Proposition~\ref{Prop:4:2:WH} below). 

Due both to the intrinsic smallness
of the quantity $\tau$ which is postulated in the weak Harnack inequality, and to its qualitative knowledge (even though, in principle, $\tau$ can be traced and made quantitative), it is quite natural to ask whether the value of $\al_1$ in Theorem~\ref{PL} is a sheer estimate linked to the approach we employ, and more refined tools could 
show that the actual value is $\al_1=1$. In a similar way, one wonders about the optimality of $\al_2$, since Lundstr\"om shows that $\al_2=\frac{p-N+M}{p-1}$ when $u$ is a sub-solution of the $p$-Laplacian, under the same restrictions on $N$, $M$, $p$ as in Theorem~\ref{PL-2}.

By relying on proper counterexamples, we show below that this is not the case. Indeed, having $\al_1\in(0,1)$ is a straightforward consequence of the wide generality of the De Giorgi classes, and, in particular, of the fact that they contain solutions to elliptic equations with bounded and measurable coefficients. An issue which we refrain from pursuing here is to trace the exact dependencies of $\al_1,\,\al_2$ on $\hat\gm,\,p,\,N,\,M$. 

\subsubsection{A First Counterexample}
The following example is due to Meyers \cite{meyers}. Let $x$ and $y$ be the coordinates in $\rr^2$, consider the equation
\begin{equation}\label{Eq:65}
\dvg({\mathbb A}\,Du)=0
\end{equation}
where
$$
{\mathbb A}=
\begin{bmatrix}
    1-\left(1-\mu^2\right) \frac{y^2}{x^2+y^2}     & \left(1-\mu^2\right) \frac{x y}{x^2+y^2}  \\
    \left(1-\mu^2\right) \frac{x y}{x^2+y^2}       & 1-\left(1-\mu^2\right) \frac{x^2}{x^2+y^2}  \\
\end{bmatrix},
$$
and $\mu$ is a fixed constant in the range $0<\mu<1$. It can be easily seen that at each point $(x, y)$ the eigenvalues of the coefficient matrix are $\mu^2$ and $1$. Thus, equation \eqref{Eq:65} is a linear elliptic equation with bounded and measurable coefficients; the function
\begin{equation}\label{count-meyers}
u(x, y)=\left(x^2+y^2\right)^{\frac{\mu-1}{2}} \cdot x
\end{equation}
is a weak solution to equation \eqref{Eq:65} which vanishes for $x=0$, and it is easy to check that $Du$ is in $L_{\loc}^p(\rr^2)$ for any $1<p<\frac{2}{1-\mu}$. Since $u$ is a weak solution to \eqref{Eq:65} in $\rr^2$, it definitely belongs to $[DG]^+_p(\rr_+^2,\hat\gm)$
for a proper $\hat\gm$, and this shows that we cannot expect $\al_1=1$ in Theorem~\ref{PL}.

We can extend the previous result to all dimensions. Let $(x, y, z)$ be a point in the $N$-space, where $z$ stands for the remaining $(N-2)$ coordinates. We extend the given solution \eqref{count-meyers} and the coefficients by defining $u(x,y,z)=u(x,y)$, and do so similarly for $\mathbb A$. Then we have
\begin{equation}\label{Eq:63}
L_{x y} u+\Delta_z u=0,
\end{equation}
and equation \eqref{Eq:63} has the same ellipticity constants as the original equation. 
Once more, the growth rate of $u$ is $\al_1<1$.
\subsubsection{A Second Counterexample}
The previous counterexample can be generalized in a different way.
When $N=3$, the function $u(x_1,x_2,x_3)=x_1\left| x\right|^\alpha$ is a solution to 
\[
\dvg({\mathbb A}Du)=0\quad\text{ weakly in }\,\,\{x_1>0\},
\]
where
\[
{\mathbb A}=
\begin{bmatrix}
1-C\frac{x_2^2+x_3^2}{x_1^2+x_2^2+x_3^2}  & C \frac{x_1 x_2}{x_1^2+x_2^2+x_3^2}  & C \frac{x_1 x_3}{x_1^2+x_2^2+x_3^2}\\
\\
C \frac{x_1 x_2}{x_1^2+x_2^2+x_3^2}  & 1-C\frac{x_1^2+x_3^2}{x_1^2+x_2^2+x_3^2} & C \frac{x_2 x_3}{x_1^2+x_2^2+x_3^2}\\
\\
 C \frac{x_1 x_3}{x_1^2+x_2^2+x_3^2}  & C \frac{x_2 x_3}{x_1^2+x_2^2+x_3^2} & 1-C\frac{x_1^2+x_2^2}{x_1^2+x_2^2+x_3^2}
\end{bmatrix},
\]
with $2 C=-\alpha(\alpha+3)$, and $\alpha\in(-1,0)$. Since $C\in(0,1)$ and the eigenvalues are \(1,\,1-C,\,1-C\), 
the matrix is uniformly elliptic. As in the previous case, it is apparent that the growth exponent is $\alpha_1<1$.
\subsubsection{A Third Counterexample}
We now take $N=4$, and consider the elliptic $p$-Laplacian with $p=N=4$, that is 
\begin{equation}\label{Eq-4-lapl}
-\dvg({\mathbb A(x)} |D u|^{2}D u)=0.
\end{equation}
We set 
\begin{equation}\label{Eq:sub-sol}
u(x_1,x_2,x_3,x_4)=\frac{(x_1^2+x_2^2)^\frac{1}{3}}{(x_1^2+x_2^2+x_3^2+x_4^2)^\frac{\alpha}{2}},
\end{equation}
where $\alpha\in (0,\frac{2}{3})$, and deal with the matrix ${\mathbb A}$
given by
\[
\begin{bmatrix}
1-\frac{C(x_2^2+x_3^2+x_4^2)}{x_1^2+x_2^2+x_3^2+x_4^2} & \frac{C x_1 x_2}{x_1^2+x_2^2+x_3^2+x_4^2}  & \frac{C x_1 x_3}{x_1^2+x_2^2+x_3^2+x_4^2} & \frac{C x_1 x_4}{x_1^2+x_2^2+x_3^2+x_4^2}\\
\\
\frac{C x_1 x_2}{x_1^2+x_2^2+x_3^2+x_4^2} & 1-\frac{C(x_1^2+x_3^2+x_4^2)}{x_1^2+x_2^2+x_3^2+x_4^2} & \frac{C x_2 x_3}{x_1^2+x_2^2+x_3^2+x_4^2} & \frac{C x_2 x_4}{x_1^2+x_2^2+x_3^2+x_4^2} \\
\\
 \frac{C x_1 x_3}{x_1^2+x_2^2+x_3^2+x_4^2}  & \frac{C x_2 x_3}{x_1^2+x_2^2+x_3^2+x_4^2} & 1-\frac{C(x_1^2+x_2^2+x_4^2)}{x_1^2+x_2^2+x_3^2+x_4^2} & \frac{C x_3 x_4}{x_1^2+x_2^2+x_3^2+x_4^2} \\
 \\
\frac{C x_1 x_4}{x_1^2+x_2^2+x_3^2+x_4^2}  & \frac{C x_2 x_4}{x_1^2+x_2^2+x_3^2+x_4^2} & \frac{C x_3 x_4}{x_1^2+x_2^2+x_3^2+x_4^2} & 1-\frac{C(x_1^2+x_2^2+x_3^2)}{x_1^2+x_2^2+x_3^2+x_4^2}
\end{bmatrix},
\]
where $C=C(\alpha)=1-\frac{(2 - 3 \alpha)^2}{8}$; we have $C\in(\frac{1}{2},1)$. This matrix is a straightforward generalization of the ones considered in the previous two examples. It is not difficult to check that the eigenvalues of the matrix are \(1,\, 1-C,\, 1-C,\, 1-C\). 

We can observe that the smallest eigenvalue $1-C$ is strictly positive, due to the upper bound on $C$ listed above. Therefore, the matrix is uniformly elliptic. We want to show that $u$ in \eqref{Eq:sub-sol} is a sub-solution to the elliptic $4$-Laplacian.
We have
\begin{align*}
\frac{\partial u}{\partial x_1}
=&\frac{\frac{1}{3}(x_1^2+x_2^2)^{-\frac{2}{3}}2x_1}{(x_1^2+x_2^2+x_3^2+x_4^2)^\frac{\alpha}{2}}-\frac{\alpha(x_1^2+x_2^2)^{\frac{1}{3}}x_1}{(x_1^2+x_2^2+x_3^2+x_4^2)^\frac{\alpha+2}{2}}\\
=&\frac{1}{(x_1^2+x_2^2+x_3^2+x_4^2)^\frac{\alpha+2}{2}}\\
&\cdot\left[\frac{2}{3}(x_1^2+x_2^2)^{-\frac{2}{3}}x_1(x_1^2+x_2^2+x_3^2+x_4^2)-\alpha(x_1^2+x_2^2)^{\frac{1}{3}}x_1\right],\\
\frac{\partial u}{\partial x_2}
=&\frac{\frac{1}{3}(x_1^2+x_2^2)^{-\frac{2}{3}}2x_2}{(x_1^2+x_2^2+x_3^2+x_4^2)^\frac{\alpha}{2}}-\frac{\alpha(x_1^2+x_2^2)^{\frac{1}{3}}x_2}{(x_1^2+x_2^2+x_3^2+x_4^2)^\frac{\alpha+2}{2}}\\
=&\frac{1}{(x_1^2+x_2^2+x_3^2+x_4^2)^\frac{\alpha+2}{2}}\\
&\cdot\left[\frac{2}{3}(x_1^2+x_2^2)^{-\frac{2}{3}}x_2(x_1^2+x_2^2+x_3^2+x_4^2)-\alpha(x_1^2+x_2^2)^{\frac{1}{3}}x_2\right],\\
\frac{\partial u}{\partial x_3}
&=\frac{-\alpha(x_1^2+x_2^2)^{\frac{1}{3}}x_3}{(x_1^2+x_2^2+x_3^2+x_4^2)^\frac{\alpha+2}{2}},\\
\frac{\partial u}{\partial x_4}
&=\frac{-\alpha(x_1^2+x_2^2)^{\frac{1}{3}}x_4}{(x_1^2+x_2^2+x_3^2+x_4^2)^\frac{\alpha+2}{2}}.
\end{align*}
Consequently, 
\begin{equation*}
\begin{aligned}
|\nabla u|^2&=\left(\frac{\partial u}{\partial x_1}\right)^2+\left(\frac{\partial u}{\partial x_2}\right)^2+\left(\frac{\partial u}{\partial x_3}\right)^2+\left(\frac{\partial u}{\partial x_4}\right)^2\\
\\
&=\frac{\left[4(x_3^2+x_4^2)+(x_1^2+x_2^2)(2-3\alpha)^2\right]}{9(x_1^2+x_2^2)^{\frac{1}{3}}(x_1^2+x_2^2+x_3^2+x_4^2)^{\alpha+1}}.
\end{aligned}
\end{equation*}
Moreover,
\begin{equation*}
\begin{aligned}
|\nabla u|^2&\left[a_{11}\frac{\partial u}{\partial x_1}+a_{12}\frac{\partial u}{\partial x_2}+a_{13}\frac{\partial u}{\partial x_3}+a_{14}\frac{\partial u}{\partial x_4}\right]\\
&=\frac{-x_1\,\left[2(-1+C)(x_3^2+x_4^2)+(x_1^2+x_2^2)(-2+3\alpha)\right]}{3(x_1^2+x_2^2+x_3^2+x_4^2)^{\frac{\alpha}{2}+1}(x_1^2+x_2^2)^\frac{2}{3}},\\
|\nabla u|^2&\left[a_{21}\frac{\partial u}{\partial x_1}+a_{22}\frac{\partial u}{\partial x_2}+a_{23}\frac{\partial u}{\partial x_3}+a_{24}\frac{\partial u}{\partial x_4}\right]\\
&=\frac{-x_2\,\left[2(-1+C)(x_3^2+x_4^2)+(x_1^2+x_2^2)(-2+3\alpha)\right]}{3(x_1^2+x_2^2+x_3^2+x_4^2)^{\frac{\alpha}{2}+1}(x_1^2+x_2^2)^\frac{2}{3}},\\
|\nabla u|^2&\left[a_{31}\frac{\partial u}{\partial x_1}+a_{32}\frac{\partial u}{\partial x_2}+a_{33}\frac{\partial u}{\partial x_3}+a_{34}\frac{\partial u}{\partial x_4}\right]
=\frac{x_3(x_1^2+x_2^2)^\frac{1}{3}(2C-3\alpha)}{3(x_1^2+x_2^2+x_3^2+x_4^2)^{\frac{\alpha}{2}+1}},\\
|\nabla u|^2&\left[a_{41}\frac{\partial u}{\partial x_1}+a_{42}\frac{\partial u}{\partial x_2}+a_{43}\frac{\partial u}{\partial x_3}+a_{44}\frac{\partial u}{\partial x_4}\right]
=\frac{x_4(x_1^2+x_2^2)^\frac{1}{3}(2C-3\alpha)}{3(x_1^2+x_2^2+x_3^2+x_4^2)^{\frac{\alpha}{2}+1}}.
\end{aligned}
\end{equation*}
Finally, by the direct calculation, we obtain
\begin{equation*}
\begin{aligned}
&\sum_{i=1}^4\frac{\partial}{\partial x_i}\left\{|\nabla u|^2\left[a_{i1}\frac{\partial u}{\partial x_1}+a_{i2}\frac{\partial u}{\partial x_2}+a_{i3}\frac{\partial u}{\partial x_3}+a_{i4}\frac{\partial u}{\partial x_4}\right]\right\}\\
=&\frac{1}{27}(x_1^2+x_2^2+x_3^2+x_4^2)^{-\frac{3\alpha}{2}-2}\left\{(2-3\alpha)^2(9\alpha^2-12\alpha+4C)(x_1^2+x_2^2)\right.\\
&+\left.\left[24\alpha(-4+3\alpha)-4C(-4-12\alpha+9\alpha^2)\right](x_3^2+x_4^2)\right\},
\end{aligned}
\end{equation*}
and since $C=1 -\frac{ (2 - 3 \alpha)^2}{8}$, we conclude
\begin{equation*}
\begin{aligned}
\sum_{i=1}^4&\frac{\partial}{\partial x_i}\left\{|\nabla u|^2\left[a_{i1}\frac{\partial u}{\partial x_1}+a_{i2}\frac{\partial u}{\partial x_2}+a_{i3}\frac{\partial u}{\partial x_3}+a_{i4}\frac{\partial u}{\partial x_4}\right]\right\}\\
&=\frac{1}{54}(x_1^2+x_2^2+x_3^2+x_4^2)^{-\frac{3\alpha}{2}-1}(2-3\alpha)^4 >0, 
\end{aligned}
\end{equation*} 
for any $\alpha\in (0,\frac{2}{3})$. Since it is easy to verify that $u\in W^{1,4}(B_R(0))$ for any $\alpha\in (0,\frac{2}{3})$ and for any $R>0$, we conclude
that indeed $u$ is a sub-solution to the equation~\ref{Eq-4-lapl} in $B_R(0)$ for any $R>0$. Hence, $u\in [DG]^+_4(\rn,\hat\gm)$ for a proper 
$\hat\gm$. Notice that here the precise value of $\hat\gm$ plays no fundamental role.

Moreover, $u$ vanishes in the $2$-dimensional hyperplane ${\mathbb H}^2=\{x_1=x_2=0\}$, and it is apparent that 
\[
\mu_+(r)\approx r^{\frac23-\alpha}.
\]
Since $\frac23-\alpha<\frac23=\frac M{N-1}$, this shows that in general we cannot expect $\alpha_2$ in Theorem~\ref{PL-2} to be equal to $\frac M{N-1}$, which is the optimal exponent growth given in \cite[Theorem~4.8]{Lindqvist1985} for the $N$-Laplacian. 

As a matter of fact, there is no need to consider the $N$-Laplacian: similar counterexamples can be built for the $p$-Laplacian for any value of $p$, provided $M$ satisfies the conditions of Theorem~\ref{PL-2}.
\subsection{Structure of the Paper}
In \S~\ref{S:notation} we introduce the notation, and above all, we consider different notions of capacity, and study their relationships. \S~\ref{S:preliminary} is devoted to a collection of known results about the main properties functions in De Giorgi's classes enjoys. In particular, we prove a logarithmic estimate; In \S~\ref{S-WH} we give a new proof of the weak Harnack inequality, which heavily relies on the structure of the sets where $u$ vanishes. Finally, \S~\ref{S:proof-p} and \ref{S:proof-N} are devoted to the proofs of the main results, respectively for $1<p<N$ and for $p=N$.
\color{black}
\subsection{Ackonwledgements}
The authors are grateful to the \emph{Erwin Schr\"odinger Institut f\"ur Mathematik und Physik} of Vienna for its kind hospitality during the Workshop ``Degenerate and Singular PDEs'' held in February 2025, where this work was completed. Both S. C. and U. G. are members of the Gruppo Nazionale per l'Analisi Matematica, la Probabilità e le loro Applicazioni (GNAMPA) of the Istituto Nazionale di Alta Matematica (INdAM); S. C. is grateful to the Departments of Mathematics of the University of Pavia, for its warm hospitality; U. G. has been partly supported by Grant 2017TEXA3H 002 “Gradient flows, Optimal Transport
and Metric Measure Structures;” Z.L. has been supported by the China Scholarship Council for 1-year study at the University of Pavia.

\section{Notation and Main Tools}\label{S:notation} 
A first simplification is the extension of $u$ to $\rn$. Indeed, for $u \in [DG]^{+}_p(\rn_+,\hat\gm)$, $u$ has an upper-semicontinuous representative 
\begin{equation*}
    u^*(x)= \lim_{\rho \downarrow 0} \sup_{B_{\rho}(x)} u = \lim_{\rho \downarrow 0} \dashint_{B_{\rho}(x)} u(y)\, dy,\quad \text{for} \quad  x \in \mathcal{L}(u, \rn_+), \quad \text{and} \quad 
0 \quad \text{otherwise},  
\end{equation*} 
where $\mathcal{L}(u,\rn_+)$ is the set of the Lebesgue points of $u$. Moreover, if $u$ satisfies \eqref{boundary-condition}, then the zero-extension
\begin{equation*}
    \tilde{u}(x)=\begin{cases} u^*(x),& x \in \rn_+,\\
    0, & x \in \rn \setminus \rn_+,    
    \end{cases}
\end{equation*} 
satisfies $\tilde{u} \in W^{1,p}_{\loc}(\rn)$ 
and \eqref{DG-boundary} with the integration extended over the full-balls. This allows us to work locally in $\rn$. The same is valid when prescribing Dirichlet conditions at the boundary of $\mathbb{H}^M$, or a more general set $E$ whose complement is $p$ fat (see Remark 1.3). In order to keep the notation simple, in the sequel we will continue to denote by $u$ its zero-extension.

\vskip0.2cm \noindent 
In the following $\omega_N$ stands for the $(N-1)$-dimensional measure of the surface of the unit ball in $\rn$.

Given a measurable set $\Om\subset\rn$, for a function $f\in L^1(\Om)$ we let
$$\dashint_{\Om}f(x)\,dx\equiv|\Om|^{-1}\int_{\Om}f(x)\,dx.$$
%
Finally, when we say that a constant $C$ depends on the \emph{data}, we mean that $C$ depends on $p$, $N$.
\subsection{Variational Capacity}\label{SS:cap1}
Let $K$ be a compact subset of $\mathbb{R}^N$ and $Q$ an open subset of $\mathbb{R}^N$ containing $K$: for $p\ge1$, the $p$-capacity of the set $K$ with respect to the set $Q$ is defined by
\begin{equation}\label{Eq:var-cap}
    \operatorname{cap}_p(K,Q)=\inf\left\{\int_{Q}|D \varphi|^p dx:\,\,\varphi\in C^\infty_0(Q),\varphi\geq 1 \text{ on } K\right\}.
\end{equation}
If $U\subset Q$ is open, we let
\[
\operatorname{cap}_p(U,Q)=\sup_{K\subset U,\,\text{compact}}\operatorname{cap}_p(K,Q),
\]
and for an arbitrary set $V\subset Q$ we define
\[
\operatorname{cap}_p(V,Q)=\inf_{V\subset U\subset Q,\,U\,\text{open}}\operatorname{cap}_p(U,Q).
\]
Let us state the main properties of $\operatorname{cap}_p$. For their proofs, we refer to \cite{Frehse}, or to \cite[Chapter~2]{HKM}.
\begin{proposition}\label{Prop:cap-proper}
Let $1<p\le N$.
\begin{enumerate}
\item $\operatorname{cap}_p(\emptyset,Q)=0$;
\item Given two sets $K_1\subset K_2\subset\subset Q$, we have $\operatorname{cap}_p(K_1,Q)\leq \operatorname{cap}_p(K_2,Q)$;
\item Given two compact sets $K_1,\, K_2\subset\subset Q$, we have 
\[
\operatorname{cap}_p(K_1\cap K_2, Q)+\operatorname{cap}_p(K_1\cup K_2,Q)\leq \operatorname{cap}_p(K_1,Q)+\operatorname{cap}_p(K_2,Q);
\]
\item For every sequence $K_j$ of compact subsets of $Q$ such that $K_1\supset\dots K_j\supset K_{j+1}\dots$, we have
\[
\operatorname{cap}_p(\cap_{j=1}^{+\infty}K_j,Q) =\lim_{j\rightarrow \infty} \operatorname{cap}_p (K_j,Q);
\]
\item For every sequence $K_j$ of subsets of $Q$ such that $K_1\subset\dots K_j\subset K_{j+1}\dots$, we have
\[
\operatorname{cap}_p(\cup_{j=1}^{+\infty}K_j,Q) =\lim_{j\rightarrow \infty} \operatorname{cap}_p (K_j,Q);
\]
\item If $V=\cup_i V_i\subset Q$, we have $\operatorname{cap}_p(V,Q)\le\sum_i \operatorname{cap}_p (V_i,Q)$;
\item Let $B_r(x_o),\, B_R(x_o)\subseteq \mathbb{R}^N$ for any $x_o\in \mathbb{R}^N$, $0<r<R$. Then, when $1<p<N$, we have
    \begin{equation*}
\operatorname{cap}_p(B_r(x_o),B_R(x_o))=\omega_N |q|^{p-1} |R^q-r^q|^{1-p},\quad q=\frac{p-N}{p-1};
    \end{equation*}
when $p=N$, we have
\begin{equation*}
    \operatorname{cap}_p(B_r(x_o),B_R(x_o))=\omega_N\left(\log\frac{R}{r}\right)^{1-N};
\end{equation*}
\item Given a compact set $K\subset Q$ and $\lm>0$, we have 
\[
\operatorname{cap}_p(\lm K,\lm Q)= \lm^{N-p}\operatorname{cap}_p(K,Q)
\] 
for any $1<p<N$;
\item Given a set $A\subset Q$, we have $|A|\le C[\operatorname{cap}_p(A,Q)]^{\frac N{N-p}}$ for any $1<p<N$ and for some constant $C$ depending only on $p$ and $N$;
\item Given a set $A\subset Q$, we have $\operatorname{cap}_p(A,Q)\le{\mathcal H}^{N-p}(A)$ for any $1<p<N$, where $\mathcal H^s$ is the $s$-dimensional Hausdorff measure; the constant $C$ depends only on $p$ and $N$.
\end{enumerate}
\end{proposition}
The last statement can be refined. Once more, we refer to \cite{Frehse} for the proofs.
\begin{proposition}
Let $1<p<N$ and $A\subset\subset Q$. If ${\mathcal H}^{N-p}(A)<\infty$, then $\operatorname{cap}_p(A,Q)=0$.
On the other hand, $1<p\le N$ and $\operatorname{cap}_p(A,Q)=0$, the ${\mathcal H}^s(A)=0$ for all $s>N-p$.
\end{proposition}
\subsection{Fat Sets}\label{SS:fat}
We come to a notion that plays a fundamental role in the following.
\begin{definition}
Given a set $E\subset\rn$, for $1<p<N$ we say that $E$ is \emph{$p$-locally uniformly fat} if there exist $r_o>0$ and $c_o>0$ such that for any $0<r<r_o$ and for every $x\in E$ we have
\begin{equation}\label{Eq:fat<N}
\frac{\operatorname{cap}_p(E\cap B_r(x),B_{2r}(x))}{\operatorname{cap}_p(B_r(x),B_{2r}(x))}\ge c_o.
\end{equation}
If $p=N$, then \eqref{Eq:fat<N} is replaced by
\begin{equation*}
\operatorname{cap}_N(E\cap B_r(x),B_{2r}(x))\ge c_o.
\end{equation*}
\end{definition}
In order to understand the role played by such a notion in our context, let us first consider the case of $E=\rn_-$. 

It is easy to see that $E$ is $p$-locally uniformly fat for any $1<p\le N$ with $r_o=+\infty$. First of all, it is apparent that it suffices to consider $x$ which belongs to the hyperplane $x_N=0$. 

If we let $B_r^-(x)\df=B_r(x)\cap\rn_-$, for any $p\in(1,N)$, by (7) and (9) of Proposition~\ref{Prop:cap-proper} we have
\begin{align*}
{\operatorname{cap}_p(E\cap B_r(x),B_{2r}(x))}&={\operatorname{cap}_p(B_r^-(x),B_{2r}(x))}\ge{C(N,p)} r^{N-p}\\
&\Rightarrow \frac{\operatorname{cap}_p(E\cap B_r(x),B_{2r}(x))}{\operatorname{cap}_p(B_r(x),B_{2r}(x))}\ge C(N,p).
\end{align*}

If $p=N$, we have
\begin{align*}
\operatorname{cap}_N(E\cap B_r(x),B_{2r}(x))=\operatorname{cap}_N(B_r^-(x),B_{2r}(x))\ge\operatorname{cap}_N(J_r(x),B_{2r}(x)),
\end{align*}
where $J_r=\{(0,0,\dots,0,x_N):\,-r<x_N<0\}$; by \cite{Lindqvist1984} we have
\[
\operatorname{cap}_N(J_r(x),B_{2r}(x))>\frac{\om_{N-2}}{\kappa_N^{N-1}}\ln 3,
\]
where $\dsty\kappa_N=\int_0^{\frac\pi2}(\sin t)^{\frac{2-N}{N-1}}\,dt$, whence the conclusion immediately follows.

Now, let us consider $E={\mathbb H}^{M}$, the $M$-dimensional hyperplane through the origin, that is
\[ 
{\mathbb H}^{M}=\{(x_1,\dots,x_M,0,0,\dots,0)\}\quad\text{ with }\, 1\le M\le N-1.
\] 
Again, we have that $E$ is $p$-locally uniformly fat with $r_o=\infty$ for any $1\le M\le N-1$, provided that $N-M<p\le N$. This is proved in \cite[Lemma~3.2]{Lund}.

Moreover, when $p=N$, it is a straightforward consequence of the following estimate proven in \cite{Anderson-1974}:
\begin{equation*}
    \operatorname{cap}_{N}(E\cap B_r(x),B_{2r}(x))
    \ge \om_{q-1}\om_{N-q-1}2^{1-q} \frac{\dsty\left[\int_1^3\frac{(t^2-1)^{q-1}}{t^q} dt\right]}
    {\dsty\left[\int_{0}^\frac{\pi}{2}(\sin t)^\frac{1+q-N}{N-1}dt\right]^{N-1}},
\end{equation*}
where $x\in{\mathbb H}^M$, and $q=1,2,\cdots,N-1$.

\begin{remark}
\upshape It is worth pointing out that the explicit value of the $2$-capacity of a disk $D_r$ of radius $r$ or of a semi-ball $B^-_r$ of the same radius with respect to $\rr^3$ are given in \cite[Chapter~II, \S~3, No.~14]{Landkof}. We have
\[
\operatorname{cap}_{2}(D_r,\rr^3)=\frac{2r}{\pi^2},\qquad \operatorname{cap}_{2}(B^-_r,\rr^3)=\frac{2r}{\pi}\left(1-\frac1{\sqrt3}\right).
\]
\end{remark}

Finally, concerning $p$-locally uniformly fat sets, we have the following result. 
\begin{proposition}[Lewis, \cite{Lewis}]\label{Prop:p-fat} 
Given $1<p\le N$, suppose $E\subset \rn$ is closed and $p$-locally uniformly fat
with constants $c_o$ and $r_o$. Then, there exist $\varep_o\in(0,1)$ and $c_1>0$ depending only on $p$, $N$, $c_o$, such that whenever $x\in E$, $0<r<r_o$, and $p-\varep_o<q<p$, the set $E$ is $q$-locally uniformly fat with constants $c_1$ and $r_o$, i.e.
\begin{equation*}
\frac{\operatorname{cap}_q(E\cap B_r(x),B_{2r}(x))}{\operatorname{cap}_q(B_r(x),B_{2r}(x))}\ge c_1.
\end{equation*}
\end{proposition}
The original Lewis' proof does not give an explicit functional dependence of $c_1$ on $c_o$. However, the proof provided in \cite{Mikko}, later extended to metric spaces in \cite{BMMS} (see also \cite{Lehr}), allows to conclude that 
\begin{equation*}
c_1=C(N,p) c_o.
\end{equation*}
\subsection{A Different Notion of Capacity for 
\texorpdfstring{\(1<p<N\)}{1<p<N}}\label{SS:cap2}
An equivalent notion of capacity can be given, assuming \(Q\equiv\rn\). Let us first consider the case where $1<p<N$. We do not go into details here and limit ourselves to underlining the link between capacity in this new framework and nonlinear potential theory; for more information, we refer the interested reader to \cite{khavin}. For simplicity, we will write $\operatorname{cap}_p(K)$ instead of $\operatorname{cap}_p(K,\rn)$. 

There exists a Radon measure $\mu_K$ supported in $K$ such that
\begin{equation}\label{Eq:4:2}
\operatorname{cap}_p(K)=\int_Kd\mu_K.
\end{equation}
Such a measure, called \emph{capacitary distribution of $K$}, can be identified with an element of $[W^{1,p}(\rn)]^*$, and generates a potential ${\mathcal U}_K$ by the formula
\begin{equation}\label{Eq:4:3}
{\mathcal U}_K(x)\equiv\int_{\rn}\int_{\rn}\left\lbrace{{\frac{d\mu_K(z)}{|y-z|^{N-1}}}}\right
\rbrace^{{\frac1{p-1}}}{\frac{dy}{|x-y|^{N-1}}}.
\end{equation}
\begin{remark}
\upshape If $p=2$, by the Riesz composition formula \eqref{Eq:4:3} implies
\begin{equation*}
{\mathcal U}_K(x)\equiv\int_{\rn}{\frac{d\mu_K(z)}{|x-z|^{N-2}}},
\end{equation*}
which is, up to a normalizing constant, the classical newtonian potential generated by $\mu_K$. 
\end{remark}
By (2.10) of \cite{khavin}, the energy associated with the potential ${\mathcal U}_K$ is
\[
{\mathcal E}(\mu_K)\equiv\|{\mathcal U}_K\|^p_{W^{1,p}(\rn)}=\int_{\rn}\left\lbrace{\int_{\rn}
{\frac{d\mu_K(z)}{|x-z|^{N-1}}}}\right\rbrace^{{\frac p{p-1}}}\,dx.
\]
The connection between these quantities is that the infimum in the definition of capacity \eqref{Eq:var-cap} is achieved when 
$\varphi\equiv{\mathcal U}_K$. Therefore,
\begin{equation}\label{Eq:4:4}
\operatorname{cap}_p(K)=\int_{\rn}\left\lbrace{\int_{\rn}
{\frac{d\mu_K(z)}{|x-z|^{N-1}}}}\right\rbrace^{{\frac p{p-1}}}\,dx.
\end{equation}
Let us now state and prove a Poincar\'e--type inequality.
\begin{proposition}\label{Prop:4:2}
Let $v\in W^{1,1}(B_R)$ and assume that $v$ vanishes on a subset $E_0$ of $B_R$. Let $K$ be any compact subset of $E_0$ and let $\mu_K$ denote its capacitary distribution. There exists a constant $\gamma$ depending only upon $N$ such that for $1<p<N$
\begin{equation}\label{Eq:4:5}\operatorname{cap}_p(K)\int_{B_R}|v|\,dx\leq\gamma R^N\int_{B_R}|Dv(x)|\left\lbrace{\int_{\rn}{\frac{d\mu_K(z)}{|x-z|^{N-1}}}}\right
\rbrace\,dx.
\end{equation}
\end{proposition}
\noi \textbf{Proof.} Without loss of generality, we may assume that $v$ is continuous and non-negative. Let $z\in K\subset E_0$ and $x\in B_R\backslash E_0$. Then
$$v(x)-v(z)=\int_0^{|x-z|}{\frac d{dr}} v(z+r\nu)\,dr\leq\int_0^{|x-z|}|Dv(y)|\,dr$$
where $\nu\equiv{\frac{x-z}{|x-z|}}$, $y=z+r\nu$.

Let us integrate this with respect to the Lebesgue measure $dx$ over all $x\in B_R\backslash E_0$ and then in $d\mu_K(z)$ over all $z\in K$. Recalling that $v$ vanishes in $K$ and taking into account \eqref{Eq:4:2} we obtain
\begin{equation}\label{Eq:4:6}
\operatorname{cap}_p(K)\int_{B_R}v\,dx\le\int_K d\mu_K(z)\int_{B_R\backslash E_0}\,dx\left\lbrace{\int_0^{|x-z|}|Dv(y)|\,dr}\right\rbrace.
\end{equation}
We estimate the last two integrals for $z\in K$ fixed, by expressing them in polar coordinates with pole at $z$. If $\tilde x\equiv\tilde x(z,\nu)$ is the polar representation of $\partial B_R$ with pole at $z$, we have 
\begin{align*}
&\int_{B_R\backslash E_0}\,dx\left\lbrace{\int_0^{|x-z|}|Dv(y)|\,dr}\right\rbrace\leq\int_{B_R}\,dx
\left\lbrace{\int_0^{|x-z|}|Dv(y)|\,dr}\right\rbrace\\
&=\int_0^{\tilde x(z,\nu)}|x-z|^{N-1}d|x-z|\int_{|\nu|=1}d\nu\left\lbrace{\int_0^{|x-z|}{\frac{|Dv(y)|}{|y-z|^{N-1}}}|y-z|^{N-1}\,d|y-z|}\right\rbrace\\
&\le\int_0^{2R}|x-z|^{N-1}d|x-z|\left\lbrace{\int_{|\nu|=1}\int_0^{|x-z|}{\frac{|Dv(y)|}{|y-z|^{N-1}}}|y-z|^{N-1}\,d|y-z|\,d\nu}\right\rbrace\\
&\le\gamma(N)R^N\int_{B_R}{\frac{|Dv(y)|}{|y-z|^{N-1}}}\,dy.
\end{align*} 
%
We substitute this in \eqref{Eq:4:6} and interchange the order of integration with the aid of Fubini's theorem to obtain \eqref{Eq:4:5}. 
\hfill\bbox
\begin{remark}
\upshape Proposition~\ref{Prop:4:2} continues to hold if $B_R$ is replaced by any convex set ${\mathcal C}$ 
\end{remark}
In the following we will employ Proposition~\ref{Prop:4:2} in a special form. Indeed, let $u\in W^{1,1}(B_{2R})$ and let $l$ and $k$ be any two numbers satisfying $l>k$. Define
\[
v=
\left\lbrace
\begin{aligned}
l-k&\quad\hbox{\rm if}\ u>l,\\
u-k&\quad\hbox{\rm if}\ k\leq u\le l,\\
0&\quad\hbox{\rm if}\ u<k,
\end{aligned}
\right.
\]
and set 
$$A_{l,R}\equiv\{x\in B_R:\ u(x)>l\},\ \ \ R>0.$$
If $A_{k,R}$ is defined analogously, from \eqref{Eq:4:5} applied to $v$ we deduce
\begin{equation}\label{Eq:4:7}
\operatorname{cap}_p(K)\int_{B_R}v\,dx\leq CR^N\int_{A_{k,\frac32R}\backslash A_{l,\frac32R}}|Du(x)|\left\lbrace{\int_{\rn}{\frac{d\mu_K(z)}{|x-z|^{N-1}}}}\right
\rbrace\,dx
\end{equation}
for every compact set $K\subset\{x\in B_R:\ u(x)<k\}$.
From \eqref{Eq:4:7} we finally obtain
\begin{equation}\label{Eq:4:8}
(l-k)\operatorname{cap}_p(K)|A_{l,R}|\leq c R^N\int_{A_{k,\frac32R}\backslash A_{l,\frac32R}}|Du(x)|\left\lbrace{\int_{\rn}{\frac{d\mu_K(z)}{|x-z|^{N-1}}}}\right
\rbrace\,dx
\end{equation}
for every compact set $K\subset\{x\in B_R:\ u(x)<k\}$.
\vskip.2truecm
\noindent As a direct further consequence of Proposition~\ref{Prop:4:2} we have the following.
\begin{corollary}
Let $v\in W^{1,p}(B_R)$, $1<p< N$, satisfy the assumptions of Proposition~\ref{Prop:4:2} and let
$$\delta(K;R)\equiv{\frac{\operatorname{cap}_p(K)}{R^{N-p}}}.$$
There exists a constant $\gamma$ depending only upon $N$ such that
\begin{equation}\label{Eq:4:9}
[\delta(K;R)]^{1/p}\dashint_{B_R}|v|\,dx\leq\gamma R\left(\dashint_{B_R}|Dv(x)|^p\,dx\right)^{1/p}.
\end{equation}
\end{corollary}
\noi \textbf{Proof.}  From \eqref{Eq:4:5} by the H\"older inequality
\begin{align*}
\operatorname{cap}_p(K)&\int_{B_R}|v|\,dx\\
&\le\gamma R^N\left(\int_{B_R}|Dv(x)|^p\,dx\right)^{1/p}\left(\int_{B_R}\left\lbrace{\int_{\rn}
\frac{d\mu_K(z)}{|y-z|^{N-1}}}\right\rbrace^{{\frac p{p-1}}}\,dy\right)^{\frac{p-1}p}.
\end{align*}
Now \eqref{Eq:4:9} follows from \eqref{Eq:4:4} and the definition of $\delta(K;R)$.
\hfill\bbox
\subsection{A Different Notion of Capacity for 
\texorpdfstring{\(p=N\)}{p=N}}\label{SS:2:4}
When \(p=N\) things are more delicate.
 Estimates \eqref{Eq:4:2}--\eqref{Eq:4:3} and \eqref{Eq:4:4} continue to hold (see \cite{khavin}, but also \cite[Chapter~2]{AdHe}), but they need to be written not in terms of the Riesz kernel
 \[
 I_1\df=\frac{\gamma_1}{|x|^{N-1}},\quad \gamma_1=\frac{\Gamma\left(\frac{N-1}2\right)}{2 \pi^{\frac {N+1}2}},
 \]
 but in terms of the Bessel kernel
 \[
 G_1\df={\mathcal F}^{-1}\left\{\frac1{(1+|\xi|^2)^{\frac12}}\right\}\equiv\frac1{(2\pi^N)}\int_{\rn}\frac{e^{ix\cdot\xi}}{(1+|\xi|^2)^{\frac12}}\,d\xi,
 \]
 or also
 \[
 G_1=\frac1{(2\pi)^{\frac N2}} |x|^{-\frac{N-2}2}
 \int_0^\infty\frac{t^{\frac N2}}{(1+t^2)^{\frac12}} J_{\frac{N-2}2}(|x|t)\,dt,
 \]
 where $J_\nu$ denotes the Bessel function of order $\nu$. It is worth pointing out that 
 \[
 G_1(x)\sim I_1(x),\qquad\text{ as }\,\,|x|\to0,
 \]
 and that for any $c\in(0,1)$
 \[
 G_1(x)=O(e^{-c|x|})\qquad\text{ as }\,\,|x|\to\infty.
 \]
 Coming to Proposition~\ref{Prop:4:2}, the estimate
 \[
 \int_{B_R\backslash E_0}\,dx\left\lbrace{\int_0^{|x-z|}|Dv(y)|\,dr}\right\rbrace\le \gamma(N)R^N\int_{B_R}{\frac{|Dv(y)|}{|y-z|^{N-1}}}\,dy
 \]
 obviously continues to hold. When we substitute it into \eqref{Eq:4:6}, take into account that $|y-z|\le2R$, and interchange the order of integration with the aid of Fubini's Theorem, we can observe that
 \[
 \begin{aligned}
 \int_{K}{\frac{d\mu_K(z)}{|y-z|^{N-1}}}&=\int_{K}e^{\frac{|y-z|}{2R}} e^{-\frac{|y-z|}{2R}}{\frac{d\mu_K(z)}{|y-z|^{N-1}}}\\
&\le\gm\int_{K}e^{-\frac{|y-z|}{2R}}{\frac{d\mu_K(z)}{|y-z|^{N-1}}}\\
&\le\gamma\int_{\rn}G_1(|x-z|)d\mu_K(z).
 \end{aligned}
 \]
Hence, \eqref{Eq:4:5} becomes
\begin{equation*}
\operatorname{cap}_N(K)\int_{B_R}|v|\,dx\leq\gamma R^N\int_{B_R}|Dv(x)|\left\lbrace{\int_{\rn} G_1(|x-z|)\,d\mu_K(z)}\right\rbrace\,dx.
\end{equation*}
Once we recall that for $p=N$
\begin{equation*}
\operatorname{cap}_N(K)=\int_{\rn}\left\lbrace\int_{\rn}
G_1(|x-z|)\,d\mu_K(z)\right\rbrace^{{\frac N{N-1}}}\,dx,
\end{equation*}
it is apparent that Corollary~2.1 remains the same.
\subsection{A Final Comment About Capacities}
Due to the equivalences among the different notions of capacity discussed in Sections~\ref{SS:cap1}, \ref{SS:cap2}, \ref{SS:2:4}, in the following we will use them in an interchangeable way; in particular, given two quantities $A$ and $B$, $A\approx B$ will mean that their quotient is bounded from above and from below by two constant that depend only on the data. 
\section{Preliminary Results}\label{S:preliminary}
\subsection{The De Giorgi Lemma}
The following lemma is well-known; for its proof we refer to \cite{ladis}. We do not state it in its full generality, we rather concentrate on a specific case tailored on the problem we are working on. 
\begin{lemma}\label{Lm:3:1}
Let $\Om\subset\rn$ be an open domain, $s\geq2$, assume $u\in[DG]_p^+(\Om,\hat\gm)$, let $R>0$ and $x_o\in\Om$ be such that $B_{2R}(x_o)\subset\Om$, and define $\mu_+=\sup_{B_{2R}(x_o)} u$. There exists a number $\theta_o\in(0,1)$, depending only upon $\hat\gamma$ such that if
\[
\left|\left\{x\in B_R(x_o):\ u(x)>\mu_+\left(1-{\frac1{2^s}}\right)\right\}\right|<\theta_o|B_R|,
\]
then 
\[
u(x)\leq\mu_+\left(1-{\frac1{2^{s+1}}}\right),\qquad\forall\,x\in B_{\frac56 R}(x_o).
\]
\end{lemma}
\begin{remark}
\upshape 
The value of $\theta_o$ forces the reduction of the oscillation in the sense that we can accept a larger (namely closer to $1$) value of $\theta_o$ if we can take $s+l$ with $l>>1$ instead of $s+1$ in Lemma~\ref{Lm:3:1}. 
For the convenience of the proof below, we provide a specific form of $\theta_o$
\[
\theta_o=\min\bigg\{\frac{5^{N+1}\omega_N }{N 6^{N+1}}, \frac{\omega_N^N}{N^N6^{4N+2N^2}{\hat\gamma}^\frac{N}{p}\beta^N}\bigg\},
\]
$\beta$ is a positive constant, identical to the one in Lemma 3.5, Chapter 2 of \cite{ladis}.
\end{remark}
\subsection{Boundedness and Weak Harnack Inequality}
A function $u \in [DG]^+_p(\rn,\hat\gm)$ is locally bounded above, and it satisfies the following local estimate.
\begin{lemma}[Thm. 1 in \cite{DiBenedetto-Trudy}] \label{bddness} 
Let $u \in [DG]^+_p(\rn,\hat\gm)$. Then, for any $\sigma \in (0,1)$, $q>0$, there exists $C >0$ depending only on the data $\{N,p\}$, and $q$ such that for all fixed $\rho>0$, 
\begin{equation*} 
    \sup_{B_{\sigma \rho}}\,u_+\, \le \bigg( \frac{\hat\gamma^{\frac{N}{p}}C }{(1-\sigma)^N} \dashint_{B_{\rho}} u_+^q\, dx \bigg)^{\frac{1}{q}}\,. 
\end{equation*}
\end{lemma} 
Another fundamental tool is the Weak Harnack inequality for function in the ``complementary" De Giorgi class.
\begin{lemma}[Theorem~2 in \cite{DiBenedetto-Trudy}]\label{lemma-WH} 
Let $v \in [DG]^-_p(\rn,\hat\gm)$  be non-negative and let $\rho>0$. Then, there exists $\tau>0$ depending only on the data $\{p,N,\hat{\gamma}\}$ with the following property:
\noindent for every $\sigma, \eta \in (0,1)$ there exists a constant $C=C(\tau, \sigma, \eta, p,N, \hat{\gamma})>0$ such that 
\begin{equation}\label{WH}
    \bigg( \dashint_{B_{\sigma \rho}} v^{\tau}\, dx \bigg)^{\frac{1}{\tau}} \leq C(\tau,\sigma,\eta,p,N,\hat\gm)  \inf_{B_{\eta \rho}} v\, .
 \end{equation}
\end{lemma}
\begin{remark}
\upshape In Lemma~\ref{bddness} there is no need to assume $u\in [DG]^+_p(\rn,\hat\gm)$; it suffices to have $v\in[DG]^+_p(\Om,\hat\gm)$, where $\Om\subset\rn$ is a bounded, open subset, such that $B_{2\rho}\subset\Om$. Analogous remark holds for Lemma~\ref{lemma-WH}.
\end{remark}
\begin{remark}\label{Rmk:3:3}
\upshape It is a matter of straightforward computations to check that if $u\in W^{1,p}_{\loc}(\Om)$ is bounded above,  $\mu_+=\sup_\Om u$, and $u\in [DG]^+_p(\Om,\hat\gm)$, then $\mu_+-u\in[DG]^-_p(\Om,\hat\gm)$.
\end{remark}
\begin{remark}
\upshape The proof of Lemma~\ref{lemma-WH} relies in a fundamental way on a delicate covering argument, originally due to Krylov \& Safonov \cite{Kry-Saf}. The natural question arises, whether in the particular framework we are dealing with now a simpler proof is actually possible, in the spirit of the original Moser's approach. The answer is in the affirmative, as we will show in \S~\ref{S-WH}
\end{remark}
\subsection{A Logarithmic Estimate for \texorpdfstring{\(1<p<N\)}{1<p<N}}
Let $\Om$ be either $\rn_-$ or $\mathbb H^M$ for $1\le M\le N-1$ and $N-M<p< N$; take $u\in[DG]^+_p(\rn,\hat\gm)$,  assume that $u\ge0$, $u\not\equiv0$, and vanishes in $\Om$. 
For simplicity, in the sequel we refer to balls centered at the origin.

For $R>0$ fixed, let 
\[
Q(R)\df=\Omega\cap \overline{B_{R}},\quad
\dl(R)\df=\frac{\operatorname{cap}_p(Q(R))}{R^{N-p}}.
\]
We have the following.
\begin{proposition}\label{pro-log-p}
Under the previous assumptions, for any $\sig\in(0,\frac{p-1}p)$ the exists $C>1$ which depends on $N$, $p$, $M$, and $\sig$ such that
\[
\int_{B_R}\left[\ln \left(\frac{\mu_+}{\mu_+-u}\right)\right]^{\sigma} dx \le \frac{C\,\hat\gamma^\frac{N+1}{p}}{[\delta(R)]^{\frac{1}{p}}}\left|B_R\right|,
\]
where $\mu_+=\sup_{B_{2R}} u$.
\end{proposition}
\noi \textbf{Proof.} $Q(R)$ is a compact set and let $\mu_{Q(R)}$ be its capacitary distribution. For simplicity of notation define
\begin{equation*}
\mu_+=\sup_{B_{2R}} u, \quad A_{s,R} \equiv\left\{x \in B_{R}: u(x)>\mu_{+}\left(1-\frac{1}{2^{s}}\right)\right\}, \quad s=4,5, \ldots,
\end{equation*}
and notice that $\dsty Q(R) \subset\left\{x \in B_{R}: u(x)<\mu_{+}\left(1-\frac{1}{2^{s}}\right)\right\}$, $s=4,5,\dots$.

Applying \eqref{Eq:4:8} with $l=\mu_{+}\left(1-\frac{1}{2^{s+1}}\right)$, $k=\mu_{+}\left(1-\frac{1}{2^{s}}\right)$, we obtain
\begin{equation*}
\begin{aligned}
\frac{\mu_+}{2^{s+1}} &\operatorname{cap}_p(Q(R))\left|A_{s+1, R}\right|\\
 \le& c R^{N} \int_{A_{s,\frac32R} \backslash A_{s+1,\frac32R}}|D u(x)|\left\{\int_{\rn} \frac{d \mu_{Q(R)}(z)}{|x-z|^{N-1}}\right\} dx\\
\le& c R^{N}\left(\int_{A_{s,\frac32R} \backslash A_{s+1,\frac32R}}\left|D\left(u-\mu_{+}\left(1-\frac{1}{2^{s}}\right)\right)_{+}\right|^{p} d x\right)^{\frac{1}{p}}\\
&\cdot\left(\int_{A_{s,\frac32R} \backslash A_{s+1,\frac32R}}\left\{\int_{\rn} \frac{d \mu_{Q(R)}(z)}{|x-z|^{N-1}}\right\}^{\frac{p}{p-1}} d x\right)^{\frac{p-1}{p}}\\
\le& c R^{N}\left(\int_{B_{\frac32R}}\left|D\left(u-\mu_{+}\left(1-\frac{1}{2^{s}}\right)\right)_{+}\right|^{p} 
dx\right)^{\frac{1}{p}}\\
&\cdot\left(\int_{A_{s,\frac32R} \backslash A_{s+1,\frac32R}}\left\{\int_{\rn} \frac{d \mu_{Q(R)}(z)}{|x-z|^{N-1}}\right\}^{\frac{p}{p-1}} d x\right)^{\frac{p-1}{p}}.
\end{aligned}
\end{equation*}
On the other hand,
\begin{equation*}
\begin{aligned}
\int_{B_{\frac32R}}\left|D\left(u-\mu_{+}\left(1-\frac{1}{2^{s}}\right)\right)_{+}\right|^{p} dx &\le \frac{\hat\gamma}{R^{p}} \int_{B_{2R}}\left|\left(u-\mu_{+}\left(1-\frac{1}{2^{s}}\right)\right)_{+}\right|^{p} dx,\\
\int_{B_{2R}}\left|\left(u-\mu_{+}\left(1-\frac{1}{2^{s}}\right)\right)_{+}\right|^{p} d x &\le c\left(\frac{\mu_+}{2^{s}}\right)^{p} R^{N} .
\end{aligned}
\end{equation*}
Collecting all the previous estimates yields
\begin{equation*}
\begin{aligned}
\operatorname{cap}_p&(Q(R))\left|A_{s+1, R}\right| \\
&\leq c{\hat\gamma}^\frac{1}{p} R^{\frac{N p+N-p}{p}} \cdot\left(\int_{A_{s,\frac32R} \backslash A_{s+1,\frac32R}}\left\{\int_{\rn} \frac{d \mu_{Q(R)}(z)}{|x-z|^{N-1}}\right\}^{\frac{p}{p-1}} d x\right)^{\frac{p-1}{p}},
\end{aligned}
\end{equation*}
and raising both sides to the $\frac{p}{p-1}$ power we obtain
\begin{align*}
\left[\operatorname{cap}_p(Q(R))\right]^{\frac{p}{p-1}}&\left|A_{s+1, R}\right|^{\frac{p}{p-1}}\\
 &\le c {\hat\gamma}^\frac{1}{p-1}R^{\frac{N p+N-p}{p-1}}
\left(\int_{A_{s,\frac32R} \backslash A_{s+1,\frac32R}}\left\{\int_{\rn} \frac{d \mu_{Q(R)}(z)}{|x-z|^{N-1}}\right\}^{\frac{p}{p-1}} d x\right).
\end{align*}
We add these inequalities for $s=4,5, \ldots, s^{*}-2$, where $s^{*}$ is a positive integer which can be chosen arbitrarily large; 
since $A_{s,\frac32R} \backslash A_{s+1,\frac32R}$ are pairwise disjoint sets, the right--hand side can be majorized with a convergent series. We have
\begin{equation}\label{Eq:5:8}
\begin{aligned}
\left(s^{*}-5\right)&\left[\operatorname{cap}_p(Q(R))\right]^{\frac{p}{p-1}}\left|A_{s^{*}-1, R}\right|^{\frac{p}{p-1}}\\
& \leq c R^{\frac{N p+N-p}{p-1}}  \sum_{s=1}^{\infty} \int_{A_{s,\frac32R} \backslash A_{s+1,\frac32R}}\left\{\int_{\rn} \frac{d \mu_{Q(R)}(z)}{|x-z|^{N-1}}\right\}^{\frac{p}{p-1}} dx,
\end{aligned}
\end{equation}
and recalling \eqref{Eq:4:4} we have
\begin{align*}
\sum_{s=1}^{\infty} \int_{A_{s,\frac32R} \backslash A_{s+1,\frac32R}}\left\{\int_{\rn} \frac{d \mu_{Q(R)}(z)}{|x-z|^{N-1}}\right\}^{\frac{p}{p-1}} d x &\leq \int_{\rn}\left\{\int_{\rn} \frac{d \mu_{Q(R)}(z)}{|x-z|^{N-1}}\right\}^{\frac{p}{p-1}} d x\\
&=\operatorname{cap}_p(Q(R)).
\end{align*}
From \eqref{Eq:5:8} we deduce
\begin{align*}
\left(s^{*}-5\right)\left[\operatorname{cap}_p(Q(R))\left|A_{s^{*}-1, R}\right|\right]^{\frac{p}{p-1}} &\le c R^{\frac{N p}{p-1}} R^{\frac{N-p}{p-1}} \operatorname{cap}_p(Q(R)),\\
\left|A_{s^{*}-1, R}\right| &\le \frac{c{\hat\gamma}^\frac{1}{p} R^{N}}{\left(s^{*}-5\right)^{\frac{p-1}{p}}}\left[\frac{\operatorname{cap}_p(Q(R))}{R^{N-p}}\right]^{-\frac{1}{p}}.
\end{align*}
%
%
%
%
%
%
%
%
%
%
%
%
%
%
%
The previous inequality can be rewritten as
\begin{equation}\label{Eq:5:12}
\left|\left[u>\mu_+\left(1-\frac{1}{2^{s}}\right)\right] \cap B_{R}\right| \leq \frac{c}{(s-4)^{\frac{p-1}{p}}[\delta(R)]^{\frac1p}}\left|B_{R}\right|,\quad s\ge5.
\end{equation}
If we define
$$
v\df=\ln \left(\frac{\mu_+}{\mu_+-u}\right),
$$
taking also into account the change of base in the logarithm, \eqref{Eq:5:12} yields
$$
\begin{aligned}
    \left|[v>s] \cap B_{R}\right|=\left|\left[u>\mu_+\left(1-\frac{1}{e^{s}}\right)\right] \cap B_{R}\right|
    &\leq \left|\left[u>\mu_+\left(1-\frac{1}{2^{s}}\right)\right] \cap B_{R}\right|\\
    &\leq \frac{c\hat\gamma^\frac{N+1}{p}}{s^{\frac{p-1}{p}}[\delta(R)]^{\frac{1}{p}}}\left|B_{R}\right|.
\end{aligned}
$$

If we now pick any $\sigma \in(0,\frac{p-1}p)$, we have
$$
\begin{aligned}
\int_{B_R} v^{\sigma} dx & =\int_{0}^{+\infty} s^{\sigma-1}\left|[v>s] \cap B_R\right| ds \\
& =\int_{0}^{5} s^{\sigma-1}\left|[v>s] \cap B_R\right| ds+\int_{5}^{+\infty} s^{\sigma-1}\left|[v>s] \cap B_R\right| ds \\
& \leq \frac{1}{\sigma}\left|B_R\right|+c\,\hat{\gamma}^\frac{N+1}{p}\left|B_R\right| \int_{5}^{\infty} \frac{s^{\sigma-1}}{(s-4)^{\frac{p-1}{p}}[\delta(R)]^{\frac{1}{p}}} ds \\
& \leq\left[\frac{1}{\sigma}+\frac{c\,\hat{\gamma}^\frac{N+1}{p}}{\left(\frac{p-1}p-\sigma\right) [\delta(R)]^{\frac{1}{p}}}\right]\left|B_R\right|=\frac{C(\sigma)\hat{\gamma}^\frac{N+1}{p}}{[\delta(R)]^{\frac{1}{p}}}\left|B_R\right|,\\
\end{aligned}
$$
and we conclude. \hfill\bbox
\subsection{A Logarithmic Estimate for \texorpdfstring{$p=N$}{p=N}}
As before, we let $\Om$ be either $\rn_-$ or $\mathbb H^M$ for $1\le M\le N-1$; take $u\in[DG]^+_p(\rn,\hat\gm)$,  assume that $u\ge0$, and vanishes in $\Om$. 
\begin{proposition}\label{pro-log-N}
Under the previous assumptions, for any $\sig\in(0,\frac{N-1}N)$ the exists $C>1$ which depends on $N$, $M$, and $\sig$ such that
\[
\int_{B_R}\left[\ln \left(\frac{\mu_+}{\mu_+-u}\right)\right]^{\sigma} dx \le \frac{C\,\hat\gamma^\frac{N+1}{N}}{[\delta(R)]^{\frac{1}{N}}}\left|B_R\right|,
\]
where $\mu_+=\sup_{B_{2R}} u$.
\end{proposition}
The proof is largely similar to the one of Proposition~\ref{pro-log-p}; indeed, it suffices to substitute
\[
\left\{\int_{\rn} \frac{d \mu_{Q(R)}(z)}{|x-z|^{N-1}}\right\}^{\frac{p}{p-1}}
\]
with
\[
\left\{\int_{\rn} G_1(|x-z|)\,d \mu_{Q(R)}(z)\right\}^{\frac{N}{N-1}},
\]
$\operatorname{cap}_p(Q(R))$ with $\operatorname{cap}_N(Q(R))$, and recall that for $p=N$ we have $\delta(R)=\operatorname{cap}_N(Q(R))$.

\section{A Weak Harnack Inequality}\label{S-WH}
\subsection{De Giorgi Classes and Sub(Super)-Harmonic 
Functions}
For any open set $\Om\subset\rn$, the so-called \emph{generalized De Gior\-gi classes} $[GDG]_p^{\pm}(\Om,\gm)$ 
are the collection of functions $u\in W^{1,p}_{\loc}(\Om)$, 
for some $p>1$, satisfying 
\begin{equation*}
\int_{B_\rho(y)}|D(u-k)_\pm|^pdx\le\frac{\gm}{(R-\rho)^p}
\Big(\frac{R}{R-\rho}\Big)^{Np}\int_{B_R(y)}|(u-k)_\pm|^pdx
\end{equation*}
for all balls $B_\rho(y)\subset B_R(y)\subset \Om$, and all 
$k\in\rr$, for a given positive constant $\gm$. An equivalent definition can be given in terms of cubes, instead of balls.

Convex, 
monotone, non-decreasing functions of 
sub-harmonic functions are sub-har\-monic. Similarly, 
concave, non-decreasing, functions of super-har\-monic 
functions are super-harmonic. Similar statements hold 
for weak, sub(super)-solutions of linear elliptic equations 
with measurable coefficients \cite{moser}.  
The next lemma establishes analogous properties 
for functions $u\in [DG]^{\pm}_p(\Om,\hat\gm)$. 
\begin{lemma}[Lemma~2.1 in \cite{DiBenedetto-Gianazza2016-bis}]\label{Lm:2:1} 
Let $\vp:\rr\to\rr$ be convex and non-decreasing, 
and let $u\in [DG]_p^+(\Om,\hat\gm)$. 
There exists a positive constant $\ov{\gm}$ 
depending only on the data, $\hat\gamma$, and 
independent of $u$, such that 
$\vp(u)\in [GDG]_p^+(\Om,\ov{\gm})$. 

Likewise let $\psi:\rr\to\rr$ be concave and non-decreasing, 
and let $u\in [DG]_p^-(\Om,\hat\gm)$. 
There exist a positive constant $\ov{\gm}$ 
depending only on the data,$\hat\gamma$, and 
independent of $u$, such that 
$\psi(u)\in [GDG]_p^-(\Om,\ov{\gm})$. 
\end{lemma}
We also have the following results.
\begin{lemma}\label{Lm:2:3}
Let $u\in [DG]_p^{+}(\Om,\hat\gm)$, assume $0\le u<1$, and take $\vp(t):=\frac1{1-t}$, $\psi(t)=\ln\frac1{1-t}$ with $t\in[0,1)$. Then for all $\tau>0$ 
there exists a constant $C_\tau$ depending only upon 
the data $\{N,p\}$ and $\tau$, such that 
\begin{align}
\sup_{B_\rho(y)}\frac1{1-u}&\le C_\tau\,\left(1+\hat\gamma\right)^\frac{(N+1)N}{p\tau}
\Big(\frac{R}{R-\rho}\Big)^{\frac{N}{\tau}} 
\left(\dashint_{B_R(y)}\left(\frac1{1-u}\right)^\tau dx\right)^{\frac1\tau},\label{Eq:2:9}
\\
\sup_{B_\rho(y)}\ln\left(\frac1{1-u}\right)&\le C_\tau \,\left(1+\hat\gamma\right)^\frac{(N+1)N}{p\tau}
\Big(\frac{R}{R-\rho}\Big)^{\frac{N}{\tau}} 
\left\{\dashint_{B_R(y)}\left[\ln\left(\frac1{1-u}\right)\right]^\tau dx\right\}^{\frac1\tau},\label{Eq:2:10}
\end{align}
for any pair of balls $B_\rho(y)\subset B_R(y)\subset \Om$.
\end{lemma}
\noi\textbf{Proof.}
Since both $\varphi$ and $\psi$ are convex, strictly increasing functions, \eqref{Eq:2:9}--\eqref{Eq:2:10} are straightforward consequences of Lemma~\ref{Lm:2:1} and of Lemma~\ref{bddness}.
\hfill\bbox
\begin{proposition}[Proposition~3.3 in \cite{DiBenedetto-Gianazza2016-bis}]\label{Prop:3:1}
Let $u\in [DG]_p^-(\Om,\hat\gamma)$ be non-negative and upper bounded
by some positive constant $M$.
Then
\begin{equation*}
\int_{B_\rho(y)}|D\ln u|^pdx\le \frac{\hat\gamma p}{(R-\rho)^p}\int_{B_R(y)} 
\ln{\frac{M}{u}}dx
\end{equation*}
for any pair of balls $B_\rho(y)\subset B_{R}(y)\subset\Om$. 
\end{proposition}
\begin{remark}
\upshape An analogous statement holds on cubes as well.
\end{remark}
\subsection{A Weak Harnack Inequality}
We have all the tools to prove a weak Harnack inequality using a Moser-like approach.
\begin{proposition}\label{Prop:4:2:WH}
Let $u\in[DG]^+_p(B_R,\hat\gm)$ for some $R>0$ and $1<p\le N$, let $\mu_+\df=\sup_{B_R} u$, assume that $u\ge0$ and vanishes either in $\rn_-\cap B_R$ or in ${\mathbb H}^N\cap B_R$ with $1\le M\le N-1$ and $N-M<p\le N$. Let $\rho\in(0,\frac R2)$; then, there exists $\tau_o\in(0,1)$, such that for any $\tau\in(0,\tau_o)$ and $s\in(0,\frac14)$ we have
\[
\dashint_{B_{2s\rho}}(\mu_+-u)^{\tau} dx\le C\, \left(1+\hat\gamma\right)^\frac{(N+1)N}{p} \inf_{B_{s\rho}}(\mu_+-u)^\tau,
\]
where $C\ge 1$ depends on the data, and $\tau$.
\end{proposition}
\noi\textbf{Proof.}
Without loss of generality, we may assume $\mu_+\equiv1$: indeed, since De Giorgi classes are homogeneous with respect to $u$, by a simple normalization such a condition is always satisfied. By Lemma~\ref{Lm:2:3}, for any $\sig>0$ we have
\[
\begin{aligned}
\int_{B_\rho}\ln\left(\frac1{1-u}\right)\,dx
&\le C_N\rho^N\sup_{B_\rho}\ln\left(\frac1{1-u}\right)\\
&\le C_N\rho^N C_\sig
\left(1+\hat\gamma\right)^\frac{N(N+1)}{p\sigma}
2^{\frac{N}{\sig}} 
\left[\dashint_{B_{2\rho}}\left[\ln\left(\frac1{1-u}\right)\right]^\sig dx\right]^{\frac1\sig},
\end{aligned}
\]
that is
\[
\left[\dashint_{B_\rho}\ln\left(\frac1{1-u}\right)\,dx\right]^\sig\le C_1 
\left(1+\hat\gamma\right)^\frac{N(N+1)}{p}
\dashint_{B_{2\rho}}\left[\ln\left(\frac1{1-u}\right)\right]^\sig dx,
\]
where $C_1$ depends on the data and $\sig$. On the other hand, by Propositions~\ref{pro-log-p}--\ref{pro-log-N}, provided we choose $\sig\in(0,\frac{p-1}p)$, we have
\[
\int_{B_{2\rho}}\left[\ln \left(\frac{1}{1-u}\right)\right]^{\sigma} dx \le \frac{C\left(1+\hat\gamma\right)^\frac{N+1}{p}}{[\delta(2\rho)]^{\frac{1}{p}}}\left|B_{2\rho}\right|,
\]
and by the argument of \S~\ref{SS:fat} there exists $c_o$ which depends only on the data, such that $\dsty \delta(2\rho)\approx\frac{\operatorname{cap}_p(Q(2\rho),B_{4\rho})}{(2\rho)^{N-p}}\ge c_o$  for any $\rho\in(0,\frac R2)$. 
Thus, for any $\sig\in(0,\frac{p-1}p)$, for any $R>0$, and for any $\rho\in(0,\frac R2)$, 
\[
\dashint_{B_\rho}\ln\left(\frac1{1-u}\right)\,dx\le C_2 \left(1+\hat\gamma\right)^\frac{(1+N)^2}{p\sigma},
\]
where $C_2$ depends on the data and $\sig$; moreover, since $u\in[DG]^+_p(B_R,\hat\gm)$, by Remark~\ref{Rmk:3:3}
we have $1-u\in [DG]^-_p(B_R,\hat\gm)$ and by Proposition~\ref{Prop:3:1}
\[
\int_{B_{s\rho}}|D  \ln(1-u)|^p dx\le\frac {C_3 \hat\gamma\left(1+\hat\gamma\right)^\frac{(1+N)^2}{p\sigma}}{(1-s)^p} \rho^{N-p},\qquad s\in(0,\frac{1}{4}).
\]
Working as in \cite[Lemma~1.2]{Trud} we conclude that there exist $C_4=4$ and 
\[
\tau_o=\frac{\ln2}{2^{N+1}C_3 {\hat\gamma}^\frac{1}{p}(1+\hat\gamma)^\frac{(1+N)^2}{p^2\sigma}}\in(0,1),
\]
such that for all $\tau\in(0,\tau_o)$ and for all $s\in(0,\frac14)$
\[
\begin{aligned}
\left(\dashint_{B_{2s\rho}}(1-u)^{-\tau} dx\right)&\left( \dashint_{B_{2s\rho}}(1-u)^{\tau} dx\right)\le C_4\\
&\quad\Leftrightarrow\quad \dashint_{B_{2s\rho}}(1-u)^{\tau} dx\le\frac {C_4}{\dsty \dashint_{B_{2s\rho}}\frac1{(1-u)^{\tau}} dx}.
\end{aligned}
\]
On the other hand, \eqref{Eq:2:9} can be rewritten as
\[
\frac1{\inf_{B_{s\rho}}(1-u)^\tau}\le {C_5\left(1+\hat\gamma\right)^\frac{(N+1)N}{p}}\,\dashint_{B_{2s\rho}}\frac1{(1-u)^\tau}dx,
\]
where $C_5$ depends on the data and $\tau$.
Combining the two previous inequalities, we conclude. 
\hfill\bbox
\section{Proof of Theorems~\ref{PL}--\ref{PL-2} for 
\texorpdfstring{\(1<p<N\)}{1<p<N}}\label{S:proof-p}
\vskip.2truecm
We give a single proof of both statements. 
%
Without loss of generality, from here on we assume that $u\in[DG]^+_p(\rn,\hat\gm)$, $u\ge0$, $u\not\equiv0$ in $\rn$, and $u$ vanishes either in $\Om\df=\rn_-$ or in $\Om\df=\mathbb H^M$ with $1\le M\le N-1$ and $N-M<p<N$. As discussed in \S~\ref{SS:fat}, in both instances the vanishing set $\Om$ is a $p$-locally uniformly fat set.

It is a matter of straightforward computations to show that any $u \in W^{1,s}\left(B_{2R}\right)$, 
$s\in\left[\frac{Np}{N+p}, p\right]$ satisfies the following Sobolev-Poincar\'e in equality
\begin{equation*}
\left(\int_{B_{\tfrac32R}} |u|^{p}\,dx\right)^{\frac1p} \le C(p, s) \frac{R^{\frac Np}}{\left[\operatorname{cap}_{s}\left(Q(\tfrac32R), B_{2R}\right)\right]^{\frac1s}}\left(\int_{B_{2R}}|Du|^{s}\, dx\right)^{\frac1s},
\end{equation*}
where $Q(\tfrac32R)=\left\{x \in B_{\tfrac32R}\cap\Om\,\text{ s.t. }\,\,u(x)=0\right\}$. See, for example \cite{ladis}, or \cite[Theorem~2.9]{Frehse}. It is apparent that $Q(\tfrac32R)\subset\Om$.

Taking into account the previous results, we conclude that
$$
\left(\int_{B_{\frac32 R}}|u|^{p}\,dx\right)^{\frac1p} \le C_{1}(N,p,q) \frac{R^{\frac Np}}{\left[\operatorname{cap}_{q}\left(Q(\frac32 R),B_{2 R}\right)\right]^{\frac1q}}\left(\int_{B_{2R}}|Du|^{q}\,dx\right)^{\frac1q},
$$
for any $p-\varepsilon_o<q<p$, where $\varep_o$ is the quantity introduced in Proposition~\ref{Prop:p-fat}. 

By the same proposition, in particular, relying both on the dependence of the $q$-fatness on the original $p$-fatness, and on the explicit expression of the $q$-capacity of a ball yields

\begin{equation}\label{Poin-fat}
\left(\int_{B_{{\tfrac32}R}}|u|^{p}\,dx\right)^{\frac1p} \le C_{2} \frac{R^{\frac Np}}{R^{\frac{N-q}{q}}}\left(\frac{1}{c_o}\right)^{\frac1q}\left(\int_{B_{2R}} |Du|^{q}\,dx\right)^{\frac1q}.
\end{equation}
Finally, since $u \in[DG]_{p}^{+}\left(\mathbb{R}^{N},\hat\gm\right)$, and $u\ge0$, by Definition~\ref{DGClass} with $k=0$ we have
$$
\int_{B_{R}}\left|Du\right|^{p}\,dx \le \frac{\hat\gamma}{(1-\tfrac23)^{p}} \frac{1}{R^{p}} \int_{B_{\tfrac32 R}} u^{p}\,dx,
$$
and we can conclude that
$$
\left(\int_{B_{R}} |Du|^{p}\,dx\right)^{\frac1p} \le C_{3}{\hat\gamma}^\frac{1}{p}  R^{\frac Np-\frac{N}q}\left(\frac{1}{c_o}\right)^{\frac1q}\left(\int_{B_{2R}}\left|Du\right|^{q}\,dx\right)^{\frac1q},
$$
even though we are not going to use this result.
Therefore, we have proved a reverse Holder inequality, much as already done by Kilpel\"ainen \& Koskela (\cite{KK}).

We let $\dsty \mu_+(r)\df=\sup_{B_r}u$, choose $\ell=\frac{7}{8} \mu_+(4 R)$, 
and define
$$
v:= 
\begin{cases}u & \text { if } 0<u \le \ell, \\ 
\ell & \text { if } u>\ell.
\end{cases}
$$
Then, \eqref{Poin-fat} yields
$$
\left(\int_{B_{\tfrac32 R}} v^{p}\,dx\right)^{\frac1p} \le C_2 \frac{R^{\frac Np}}{R^{\frac{N-q}{q}}}\left(\frac{1}{c_o}\right)^{\frac1q}\left(\int_{B_{2R}}|Dv|^{q}\,dx\right)^{\frac1q},
$$
and also
$$
\begin{aligned}
\ell\left|A_{\ell,\tfrac{3}{2}R}\right|^{\frac1p} &\le C_{2} R^{1+\frac Np-\frac{N}{q}}\left(\int_{B_{2R}} |Du|^{q}\,dx\right)^{\frac1q}\\
\ell^{q}\left|A_{\ell,\tfrac{3}{2}R}\right|^{\frac qp} &\le C_2^q R^{q+\frac qp N-N} \int_{B_{2R}}|Du|^{q}\,dx \\
& \le C_2^q R^{q+\frac qp N-N} \int_{B_{2R}} \frac{\left|Du\right|^{q}}{(\mu_+(4R)-u)^{(1+\alpha)q}}(\mu_+(4R)-u)^{(1+\alpha)q}\, dx,
\end{aligned}
$$
where $\alpha>0$ is to be chosen. Since
$$
D \frac{1}{(\mu_+(4R)-u)^{\alpha}}=\alpha \frac{1}{(\mu_+(4R)-u)^{1+\alpha}} Du,
$$
we can rewrite
$$
\begin{aligned}
{}[\mu_+(4R)]^{q}&\left|A_{\ell,\tfrac{3}{2} R}\right|^{\frac qp}\\
&\le C_{4} R^{q+\frac {qN}p-N} \int_{B_{2R}}\left|D \frac{1}{(\mu_+(4R)-u)^{\alpha}}\right|^{q}(\mu_+(4R)-u)^{(1+\alpha)q}\,dx,
\end{aligned}
$$
where $C_4=C_4\left(\frac78,c_o,N,p,\alpha\right)$. We stipulate that
$$
\left|A_{\ell,\frac{3}{2} R}\right| \ge \theta_o\left|B_{\frac{3}{2} R}\right|,
$$
where $\theta_o$ is the quantity introduced in Lemma~\ref{Lm:3:1}; we will later consider what happens when such a stipulation does not hold. Then, for a constant $C_5$ which depends on the original choice of $\ell$, and on $c_o$, $N$, $p$, $\alpha$, we have
$$
\begin{aligned}
{}[\mu_+&(4R)]^{q} R^{\frac{qN}p}\\
\le& C_{5} {\hat\gamma}^{\frac{Nq}{p^2}} R^{q+\frac{qN}p-N} \int_{B_{2R}}\left|D\frac{(\mu_+(4R))^{\alpha}}{(\mu_+(4R)-u)^{\alpha}}\right|^{q} \frac{(\mu_+(4R)-u)^{(1+\alpha)q}}{(\mu_+(4R))^{\alpha q}}\, dx \\
\le& \frac{C_{5} {\hat\gamma}^{\frac{Nq}{p^2}} R^{q+\frac{qN}p-N}}{(\mu(4R)_+)^{\alpha q}}\left(\int_{B_{2R}}\left|D \frac{(\mu_+(4R))^{\alpha}}{(\mu_+(4R)-u)^{\alpha}}\right|^{p}\,dx\right)^{\frac qp}\\
&\cdot\left(\int_{B_{2R}}(\mu_+(4R)-u)^{(1+\alpha)\frac{qp}{p-q}}\,dx\right)^{\frac{p-q}{p}},
\end{aligned}
$$
that is, taking Lemma~\ref{Lm:2:1} into account, and choosing $\sig_o\in(0,\tau_o)$ where $\tau_o$ is the quantity stipulated in Proposition~\ref{Prop:4:2:WH}
$$
\begin{aligned}
{}[\mu_+(4R)]^{(1+\alpha)q} \le & C_{6} \left(1+\hat\gamma\right)^\frac{Nq+pq(N+1)}{p^2}{}\frac1{R^{N}}\left(\int_{B_{3 R}}\left|\frac{\mu_+(4R)}{\mu_+(4R)-u}\right|^{\alpha p}\,dx\right)^{\frac qp}\\
\cdot & \left(\int_{B_{3 R}}[\mu_+(4R)-u]^{\sigma_o}\,dx\right)^{\frac{p-q}{p}}[\mu_+(4R)]^{(1+\alpha)q-\sigma_o \frac{p-q}{p}};
\end{aligned}
$$
by Proposition~\ref{Prop:4:2:WH} this yields
$$
\begin{aligned}
{}[\mu_+(4R)]^{\sigma_o\frac{p-q}{p}} \le& C_{7} \left(1+\hat\gamma\right)^\frac{(p-q)N^2+pN+pq(N+1)}{p^2} \frac{[\mu_+(4R)]^{\alpha q}}{[\mu_+(4R)-\mu(3R)]^{\alpha q}}\\
&\cdot[\mu_+(4R)-\mu(3R)]^{\sigma_o\frac{p-q}{p}} \\
[\mu_+(4R)]^{\sigma_o \frac{p-q}{p}-\alpha q} \le& C_{7}
\left(1+\hat\gamma\right)^\frac{(p-q)N^2+pN+pq(N+1)}{p^2} [\mu_+(4R)-\mu_+(3R)]^{\sigma_o \frac{p-q}{p}-\alpha q},
\end{aligned}
$$
and we may finally choose $\alpha$ small enough, so that
$$
\sigma_o \frac{p-q}{p}-\alpha q>0\,\,\Rightarrow\,\,\alpha<\sigma_o \frac{p-q}{pq}\,\,\Rightarrow\,\,
\alpha<\sigma_o\left(\frac{1}{q}-\frac{1}{p}\right).
$$
Correspondingly, we obtain
$$
\mu_+(4R) \le C_{8}
\left(1+\hat\gamma\right)^\frac{(p-q)N^2+pN+pq(N+1)}{p(\sigma_o(p-q)-\alpha pq)} [\mu_+(4R)-\mu_+(3R)],
$$
and also
$$
\mu_+(3R) \le\left(1-\frac{1}{C_{9}}\right) \mu_+(4R),
$$
where $C_{9}=C_{8}\left(1+\hat\gamma\right)^\frac{(p-q)N^2+pN+pq(N+1)}{p(\sigma_o(p-q)-\alpha pq)}$, $C_{8}$ depends on the data, $q$, $\alpha$, $\sigma_o$, $c_o$.

On the other hand, if the previous assumption is violated, that is if $\left|A_{\ell,\frac32 R}\right|<\theta_{0}\left|B_{\frac{3}{2} R}\right|$, then by Lemma~\ref{Lm:3:1}, we have
$$
\mu_+(R) \le \frac{15}{16} \mu_+(2R).
$$
Since $\dsty \mu_+(3R) \ge \mu_+(R)$, and $\dsty\mu_+(2R) \le \mu_+(4R)$,
if we let
$$
\eta=\max \left\{\frac{15}{16}, 1-\frac{1}{C_{9}}\right\},
$$
we conclude that
$$
\mu_+(R) \le\eta\, \mu_+(4R).
$$
The constants $C_i$ depend on the set where $u$ vanishes, but the structure of the proof is independent of the size of $Q(R)$.

Let us define the sequence of radii
\[ 
R_k= 4^k, \qquad k \in \nn\,,
\]
Then, for any $r\ge 4$, there exists $k\in \mathbb{N}$, such that
\[
R_k\leq r\leq R_{k+1}=4^{k+1},
\]
which means $k\ge \frac{\ln r}{\ln 4}-1$. Therefore, we derive
$$
\begin{aligned}
    \mu_+(r)\ge \mu_+(R_{k})&\ge \frac{1}{\eta}\mu_+(R_{k-1})\\
    &\vdots\\
    &\ge \frac{1}{\eta^k} \mu_+(R_0)\ge \frac{1}{\eta^{\frac{\ln r}{\ln 4}-1}} \mu_+(R_0)\\
    &= \mu_+(R_0) \eta \exp{\left(\ln r \cdot \frac{\ln\frac{1}{\eta}}{\ln 4}\right)}\\
    &= \mu_+(R_0) \eta\, r^\frac{\ln\frac{1}{\eta}}{\ln 4}=\frac{\mu_+(R_0)}{R_0^\frac{\ln\frac{1}{\eta}}{\ln 4}}\eta\, r^\frac{\ln\frac{1}{\eta}}{\ln 4}.
\end{aligned}
$$
Consequently, we obtain
$$
\frac{\mu_+(r)}{r^\frac{\ln\frac{1}{\eta}}{\ln 4}}\ge \frac{\mu_+(R_0)}{R_0^\frac{\ln\frac{1}{\eta}}{\ln 4}}\eta,
$$
and we finish the proof, once we assume $\alpha_i=\frac{\ln\frac{1}{\eta}}{\ln 4}$. 
\hfill\bbox
\begin{remark}
\upshape When $u$ vanishes in $\rn_-$, a simpler proof can be given. 
For any fixed $R>0$, let us consider the function 
\[
0 \leq (\mu_+(R)-u)\in [DG]^-_p(B_R,\hat\gm)\,.
\]
Applying Proposition \ref{Prop:4:2:WH} and using that $u$ vanishes in $B_{2s \rho}^-= B_{2s\rho} \cap \rn_-$, we obtain
  \begin{equation*}
      \begin{aligned}
          \frac{\mu_+(R)}{2^{1/\tau}} = \bigg( \frac{\mu_+(R)^{\tau}}{|B_{2s\rho}|} \bigg( |B_{2s\rho}|/2\bigg) \bigg)^{1/\tau}&= \bigg( \frac{1}{|B_{2s\rho}|} \int_{B_{2s\rho}^-} (\mu_+(R))^{\tau}\, dx \bigg)^{1/\tau}\\
          & \leq \bigg( \dashint_{B_{2s\rho}} (\mu_+(R)-u(x))^{\tau}\, dx\bigg)^{1/\tau}\\
          &\leq C \left(1+\hat\gamma\right)^\frac{(N+1)N}{\tau p} \, \, \inf_{B_{s\rho}} \bigg(\mu_+(R)-u(x)\bigg)\\
          &=C \left(1+\hat\gamma\right)^\frac{(N+1)N}{p\tau}\bigg(\mu_+(R)-\mu_+(s\rho)\bigg)\, ,
      \end{aligned}
  \end{equation*} 
  where the parameters $s,~\rho$ and $\tau$ are the same as those in Proposition \ref{Prop:4:2:WH}, and $C$ depends on the data and $\tau$.
  Hence, we have
  \[ \frac{\mu_+(R)}{2^{1/\tau}} \leq C \left(1+\hat\gamma\right)^\frac{(N+1)N}{p\tau}\, \bigg( \mu_+(R)-\mu_+(s\rho)\bigg) ,\]
which is 
\[ \eta\,\mu_+(R) \ge \mu_+(s\rho), \qquad \text{with} \qquad \eta= \bigg( 1-\frac{1}{2^{1/\tau} C\left(1+\hat\gamma\right)^\frac{(N+1)N}{p\tau}}\bigg)\,.\]
Then, we set $R_k=\left(\frac{R}{s\rho}\right)^k,~k\in \mathbb{N}$. Through a similar iterative process as in the final part of the proof above, we derive
$$
\frac{\mu_+(r)}{r^{\alpha_1}}\ge \frac{\mu_+(R_0)}{R_0^{\alpha_1}}\eta, \quad\text{where}\quad\alpha_1=\frac{\ln\frac{1}{\eta}}{\ln \frac{R}{s\rho}}.
$$
\end{remark}
\section{Proof of Theorems~\ref{PL}--\ref{PL-2} for \texorpdfstring{\(p=N\)}{p=N}}\label{S:proof-N}
When we deal with $p=N$, the gist of the argument remains essentially the same as discussed in the previous section. We only need to substitute $p$ with $N$ and take $q\in (N-\varep_o,N)$, instead of $q\in(p-\varep_o,p)$.

\section{Final Remarks}
The weak Harnack inequality of Proposition~\ref{Prop:4:2:WH} does not depend on the special choices $\Omega=\rn_-$ or $\Omega= \mathbb{H}^M$, and it can be formulated for any $0\leq u \in [DG]^+_p(\rn, \gamma)$ that vanishes in of a $p$-fat set $\Omega$.

Moreover, the proof of Theorems \ref{PL}-\ref{PL-2} relies on three different features:
\begin{itemize}
\item the zero-extension of $u$ in $\Omega$;
\item the weak Harnack inequality for $p$-fat sets; 
\item the possibility to repeat the procedure for arbitrarily big radii, i.e. for $\Omega$ unbounded.
\end{itemize}
Hence, our result is valid also for the other examples given in \cite[Chapter~I, \S~6]{Landis}, i.e. the slab $\text{S}$, where
\[
\text{S}:= \{(x_1, x_2, \dots, x_N): |x_N|<h \},
\] 
the double circular cone, the one-edged cone, and many others (we refrain from giving complete definitions of these objects, and we refer the interested reader to \cite{Landis}).

Furthermore, the weak Harnack inequality for $p$-fat sets may be useful in the future also to attack properties of singular solutions, such as in  \cite[Theorem~7.40]{HKM}. Here a similar statement seems reasonable to be valid with $\alpha<0$.

\end{document}